\documentstyle[11pt]{article}

\newfont{\Bb}{msbm10 scaled\magstephalf}

\begin{document}
\noindent
   \begin{center}
  {\LARGE The Weinstein conjecture in the uniruled manifolds}
  \end{center}
\noindent
  \begin{center}
    {\large  Guangcun Lu}\footnote{Partially supported by the NNSF
   19971045  of China.}\\[5pt]
       Nankai Institute of Mathematics, Nankai University\\
               Tianjin 300071, P. R. China\\[5pt]  \vspace{-2mm}
               (E-mail: gclu@nankai.edu.cn)\\[3pt]
              \end{center}

\begin{abstract} In this note we prove the Weinstein conjecture for a class of
symplectic manifolds including the uniruled manifolds based on Liu-Tian's result.
 \end{abstract}



Since 1978 A. Weinstein proposed his famous conjecture that {\it every hypersurface
of contact type in the symplectic manifolds carries a closed characteristic}([We]),
many results were obtained (cf.[C][FHV][H][HV1][HV2][LiuT][Lu1][Lu2][V1][V2][V3])
after C.Viterbo first proved it in $(\mbox{\Bb R}^{2n},\omega_0)$ in 1986([V1]).
Not long ago Gang Liu and Gang Tian  established a deep relation between
this conjecture and the Gromov-Witten invariants and got several general
results as corollaries([LiuT]). 

Assume $S$ to be a hypersurface of contact type in a closed connected symplectic
manifold $(V,\omega)$ separating $V$  in the sense of [LiuT], i.e. there
exist submanifolds $V_+$ and $V_-$ with common boundary $S$ such that
$V=V_+\cup V_-$ and $S=V_+\cap V_-$, then the following 
result holds.

\noindent{\bf Theorem 1}([LiuT])\quad{\it If there exist $A\in H_2(V;\mbox{\Bb Z})$
and $\alpha_+,\, \alpha_-\in H_\ast(V;\mbox{\Bb Q})$, such that 
\begin{description}
\item[(i)] $supp(\alpha_+)\hookrightarrow int(V_+)$ and
$supp(\alpha_-)\hookrightarrow int(V_-)$, 
\item[(ii)] the Gromov-Witten invariant 
$\Psi_{A, g, m+2}(C;\alpha_-,\alpha_+,\beta_1,\cdots,\beta_m)$ $\ne 0$
for some $\beta_1,\cdots,\beta_m\in H_\ast(V;\mbox{\Bb Q})$ and
$C\in  H_\ast(\overline{\cal M}_{g, m+2};\mbox{\Bb Q})$,
\end{description}
then $S$ carries at least one closed characteristic.}

Recall that for a given $A\in H_2(V;\mbox{\Bb Z})$ the Gromov-Witten invariant
of genus $g$ and with m+2 marked points is a homomorphism
$$\Psi_{A, g, m+2}: H_\ast(\overline{\cal M}_{g, m+2};\mbox{\Bb Q})\times
H_\ast(V;\mbox{\Bb Q})^{m+2}\to \mbox{\Bb Q},$$
(see [FO][LiT][R][Si]). Though one so far  does not yet know whether the GW invariants
defined in the four papers agree or not, we believe that they have the same
vanishing or nonvanishing properties, i.e., for any given classes
$C\in H_\ast(\overline{\cal M}_{g, m+2};\mbox{\Bb Q})$ and
$\beta_1,\cdots,\beta_{m+2}\in H_\ast(V;\mbox{\Bb Q})$ one of this four versions
vanishes on $(C;\beta_1,\cdots,\beta_{m+2})$ if and only if any other three ones
vanish on them. In addition, the version of [R] is actually a homomorphism from
$H_\ast(\overline{\cal M}_{g, m+2};\mbox{\Bb R})\times
H_\ast(V;\mbox{\Bb R})^{m+2}$ to $\mbox{\Bb R}$. However, using the 
facts that $H_\ast(M;\mbox{\Bb Q})$ is dense $H_\ast(M;\mbox{\Bb R})$
for $M=V, \overline{\cal M}_{g,k}$ and  $\Psi_{A, g, m+2}$ is always
a homomorphism one can naturally extend the other three versions to the
homomorphisms from $H_\ast(\overline{\cal M}_{g, m+2};\mbox{\Bb R})\times
H_\ast(V;\mbox{\Bb R})^{m+2}$ to $\mbox{\Bb R}$. Below we always mean the extended
versions when they can not clearly explained in the original versions.
Our main result is 

\noindent{\bf Theorem 2}\quad{\it For a connected closed symplectic manifold 
$(V,\omega)$, if there exist  $A\in H_2(V;\mbox{\Bb Z})$,
$C\in  H_\ast(\overline{\cal M}_{g, m+2};\mbox{\Bb Q})$ 
 and $\beta_1,\cdots, \beta_{m+1}\in H_\ast(V;\mbox{\Bb Q})$ such that
 $$\Psi_{A, g, m+2}(C;[pt],\beta_1,\cdots,\beta_{m+1})\ne 0$$ 
for $(g, m)\ne (0, 0)$ and the single  point class $[pt]$,  
then every hypersurface of contact type $S$ in the symplectic manifold $V$
separating $V$ carries a closed characteristic. Specially, if $g=0$ we can also 
guarantee that $S$ carries one contractible ( in $V$) closed characteristic.}

In case $g=0$ it is not difficult to prove that Proposition 2.5(5) and Proposition 2.6
in [RT] still hold for any closed symplectic manifold $(V,\omega)$
with the method of [R]. That is, 
\begin{description}
\item[(i)] $\Psi_{0, 0, k}([pt]; \alpha_1,\cdots,\alpha_k)=
\alpha_1\cap\cdots\cap\alpha_k$( the intersection number);
\item[(ii)] for the product manifold $(V,\omega)=(V_1\times V_2,
\omega_1\oplus\omega_2)$ of any two closed symplectic manifolds $(V_1,\omega_1)$ and
$(V_2,\omega_2)$ it holds that
$$\Psi^V_{A_1\otimes A_2, 0, k}([pt];\alpha_1\otimes\beta_1,\cdots,
\alpha_k\otimes\beta_k)=
\Psi^{V_1}_{A_1, 0, k}([pt];\alpha_1,\cdots,\alpha_k)
\Psi^{V_2}_{A_2, 0, k}([pt];\beta_1,\cdots,\beta_k).$$
\end{description}
Thus if  $\Psi^{V_2}_{A_2, 0, m+1}([pt]; [pt],\beta_1,\cdots,\beta_m)\ne 0$ we get
$$\Psi^V_{A_1\otimes A_2, 0, m+1}([pt]; [pt], \alpha_1\otimes\beta_1,\cdots,
\alpha_m\otimes\beta_m)\ne 0$$
for $A_1=0$ and $\alpha_1=\cdots =\alpha_m=[V_1]$. This leads to

\noindent{\bf Corollary 3}\quad{\it Weinstein conjecture holds in the product 
symplectic manifolds  of any closed symplectic manifold and a symplectic manifold
satisfying the condition of Theorem 2 for $g=0$.}

Recall that a smooth Kahler manifold $(M,\omega)$ is called {\it uniruled }
if it can be covered by rational curves. Y. Miyaoka and S. Mori showed that
a smooth complex projective manifold $X$ is uniruled if and only if there exists
a non-empty open subset $U\subset X$ such that for every $x\in U$ there
is an irreducible curve $C$ with $(K_X, C)<0$ through $x$([MiMo]). Specially,
any {\it Fano} manifold is uniruled([Ko]). The complex projective spaces,
the complete intersections in it, the Grassmann manifolds and more general flag manifold
are the important examples 
of the Fano manifolds. In [R, Prop. 4.9] it was proved that
if a smooth Kahler manifold $M$ is symplectic deformation equivalient
to uniruled manifold, $M$ is uniruled. Actually, as mentioned there,
Kollar showed that on the 
uniruled manifold $(M,\omega)$ there exists a class $A\in H_2(V;\mbox{\Bb Z})$
such that
$$\Phi_{A, 0, 3}([pt];[pt],\beta_1,\beta_2)\ne 0\leqno(1)$$
for some classes $\beta_1$ and $\beta_2$(see [R] for more general case).
Combing these with Corollary 3 we get

\noindent{\bf Corollary 4}\quad{\it Every hypersurface $S$ of contact type 
in the uniruled manifold $V$ or the product of any closed symplectic
manifold and an uniruled manifold carries one contractible ( in $V$) closed 
characteristic.}

The ideas of proof are combing Liu-Tian's Theorem 1 above, the properties
of the Gromov-Witten invariants and  Viterbro's trick of [V4].

\noindent{\bf Proof of Theorem 2}\quad Under the assumptions of Theorem 2,
the reduction formula of the Gromov-Witten invariants([R, Prop. C]) implies that
$$\Psi_{A, g, m+3}(\pi^\ast(C);[pt], PD([\omega]),\beta_1,\cdots,\beta_{m+1})=
\omega(A)\cdot\Psi_{A, g, m+2}(C;[pt],\beta_1,\cdots,\beta_{m+1})\ne 0\leqno(2)$$ 
since $A$ contains the nontrivial pseudoholomorphic curves. To use Theorem 1
we need to show that there exists a homology class $\gamma\in H_2(V;\mbox{\Bb R})$
with support $supp(\gamma)\hookrightarrow int(V_+)$( or $int(V_-)$) such that 
$$\Psi_{A, g, m+3}(\pi^\ast(C);[pt], \gamma,\beta_1,\cdots,\beta_{m+1})\ne 0.
\leqno(3)$$
To this goal we note  $S$ to be a hypersurface of contact 
type, and thus there exists a Liouville vector field $X$ defined in a 
neighborhood $U$ of $S$, which is transverse to $S$. The flow of $X$ define 
a diffeomorphism $\Phi$
from $S\times (-3\epsilon, 3\epsilon )$ onto an open neighborhood of $S$ in $U$
for some $\epsilon>0$. Here we may assume 
$\Phi(S\times (-3\epsilon, 0])\subset V_+$ and
$\Phi(S\times [0, 3\epsilon))\subset V_-$.
 For any $0<\delta<3\epsilon$ let us denote by
 $U_\delta:=\Phi(S\times [-\delta,\delta])$. We also denote by $\alpha=i_X\omega$, 
then $d\alpha=\omega$ on $U$. Choose a smooth function $f: V\to \mbox{\Bb R}$ such 
that $f|_{U_\epsilon}\equiv 1$ and vanishes outside $U_{2\epsilon}$. Define
$\beta:=f\alpha$. This is a smooth $1$-form on $V$, and $d\beta=\omega$ 
on $U_\epsilon$.
Denote by $\widehat\omega=\omega-d\beta$. Then 
$\widehat\omega|_{U_\epsilon}\equiv 0$ and thus cohomology classes
$[\omega]=[\widehat\omega]$ is in
$H^2(V, U_\epsilon)$.
Now from the naturality of Poincare-Lefschetz duality ([p.296, Sp]):
$H_{2n-2}(V-U_\epsilon)\cong H^{2}(V, U_\epsilon)$
  it follows that we  can choose a cycle
representive $\widehat\gamma$ of $\gamma:=PD([\omega])$  with support
$supp(\widehat\gamma)\hookrightarrow int(V-U_\epsilon)$.
Notice that $V-U_\epsilon\subset V-S=int(V_+)\cup int(V_-)$ and
$int(V_+)\cap int(V_-)=\emptyset$. We can denote by $\widehat\gamma_+$ and 
$\widehat\gamma_-$ the union of connected components of $\widehat\gamma$ lying
$int(V_+)$ and $int(V_-)$ respectively. Then the homology classes determined 
by them in $H_\ast(V,\mbox{\Bb R})$ satisfy: $[\widehat\gamma_+] +[\widehat\gamma_-]=\gamma$.
Thus $[\widehat\gamma_+]$ and $[\widehat\gamma_-]$ have at least one nonzero class. 
By the property of the Gromov-Witten invariants we get
\begin{eqnarray*}
(4)\hspace{10mm}\Psi_{A, g, m+3}(\pi^\ast(C);[pt], \gamma,\beta_1,\cdots,\beta_{m+1})
&=&\Psi_{A, g, m+3}(\pi^\ast(C);[pt], [\widehat\gamma_+],\beta_1,\cdots,\beta_{m+1})\\
&+&\Psi_{A, g, m+3}(\pi^\ast(C);[pt], [\widehat\gamma_-],\beta_1,\cdots,\beta_{m+1})
\ne 0.
\end{eqnarray*}
Hence the right side of (4) has at least one nonzero term. Without loss of
generality we assume that
$$\Psi_{A, g, m+3}(\pi^\ast(C);[pt], [\widehat\gamma_+],\beta_1,\cdots,\beta_{m+1})\ne 0.$$ 
Then Theorem 1 directly leads to the first conclusion.

The second claim is easily obtained by carefully checking the arguments in 
[LiuT].\hfill$\Box$\vspace{4mm}

\noindent{\bf Remark 5}\hspace{2mm}Actually we believe that Theorem 1 still holds
provided the hypersurface $S$ of contact type therein is replaced by the stable
hypersurface in the sense of [HV2]. Hence the hypersurface $S$ of contact type in
our results above may be replaced by the stable hypersurface for which the symplectic
form is exact in some open neighborhood of it.

\noindent{\bf Remark 6}\hspace{2mm}In [B] it was proved that the system of
Gromov-Witten invariants of the product of two varieties is equal to the
tensor product of the systems of Gromov-Witten invariants of the two factors.
Using the methods developed in [FO][LiT][R][Si] we believe that one can still
prove these product formula of Gromov-Witten invariants with any genus
for any product of two closed symplectic manifolds, and thus  Corollary 3 also holds
for any genus $g$.

{\bf Acknowledgements}\hspace{2mm}
I would like to express my hearty thanks to Professor Claude Viterbo for his many 
valuable discussions and sending his recent preprint [V4] to me. 
I wish to acknowledge Professor Janos Kollar for his telling me some
properties on the uniruled manifolds. He also thanks the referee for his good
suggestions to improve the original version.

\end{document}